\documentclass[11pt, leqno]{article}
\usepackage{amsmath,amsfonts,amssymb,wasysym,mathrsfs}
\usepackage[latin1]{inputenc}
\usepackage{shapepar}
\usepackage{graphicx}
\usepackage{hyperref}
\usepackage{ifpdf} 
\usepackage{color}
\newcommand{\R}{{\mathbb R}}
\newcommand{\vp}{\varphi}
\newcommand{\re}{{\mathbb R}}
\newcommand{\ren}{{\mathbb R}^N}

\newcommand{\be}[1]{\begin{equation}\label{#1}}
\newcommand{\ee}{\end{equation}}

\newcommand{\prf}{\par\smallskip\noindent{\sl Proof. \/}}
\newcommand{\finprf}{\unskip\null\hfill$\;\square$\vskip 0.3cm}

\newtheorem{theorem}{Theorem}[section]
\newtheorem{lemma}{Lemma}[section]
\newtheorem{corollary}[theorem]{Corollary}

\newcommand{\ve}{\varepsilon}

\numberwithin{equation}{section}
\newcommand{\nc}{\normalcolor}

\def\qed{\,\unskip\kern 6pt \penalty 500
\raise -2pt\hbox{\vrule \vbox to8pt{\hrule width 6pt
\vfill\hrule}\vrule}\par}
\definecolor{darkblue}{rgb}{0.05, .05, .65}
\definecolor{darkgreen}{rgb}{0.1, .65, .1}
\definecolor{darkred}{rgb}{0.8,0,0}

\newcommand{\nlc}{\normalcolor}
\begin{document}
\title{\textbf{ The very singular solution for the Anisotropic Fast Diffusion Equation and its consequences}\\[7mm]}

\author{
\Large Juan Luis V\'azquez \footnote{Departamento de Matem\'aticas, Universidad Aut\'onoma de Madrid, 28049 Madrid, Spain. \newline
E-mail: {\tt juanluis.vazquez@uam.es}}
}

\date{\today} 

\maketitle

\begin{abstract}
 We construct  the Very Singular Solution (VSS) for the Anisotropic Fast Diffusion Equation (AFDE) in the suitable good exponent range of  fast diffusion. VSS is a solution that, starting  from an infinite mass located at one point as initial datum, evolves  according to the corresponding equation  as an admissible solution away from the singularity. It is expected to represent important properties of the fundamental solutions when the initial mass is very big.  We work in the whole space $\ren$.

In this setting we show that  the diffusion process distributes  mass from the the initial infinite singularity  along the different space directions: up to constant factors, there is a simple partition formula for the anisotropic mass expansion, given approximately as the minimum of separate 1-D VSS solutions. This striking fact is a consequence of the improved scaling properties of the special solution, and it has strong consequences.

If we consider the family of Fundamental Solutions for different masses, we prove that they share the same universal tail behaviour as the VSS. Actually, their tail is asymptotically convergent to the unique VSS tail, so there is a VSS partition formula for their profile at spatial infinity.

 With the help of this analysis we study the behaviour of the class of nonnegative finite-mass solutions of the Anisotropic FDE, and  prove the  Global Harnack Principle (GHP)  and the Asymptotic Convergence in Relative Error (ACRE) under a natural assumption on the decay of the initial tail.

\end{abstract}

\vskip .5cm

\setcounter{page}{1}



\medskip


\

\newpage

\numberwithin{equation}{section}


\section{Introduction}\label{sec.intro}
The aim of this paper is to contribute to the mathematical theory of nonlinear evolution equations that combine a mechanism of nonlinear diffusion with the occurrence of strong anisotropy. A typical model of the equations of nonlinear diffusion is the Porous Medium Equation
\begin{equation}\label{PME}
u_t=\Delta u^m \tag{PME}
\end{equation}
posed in a domain in $\ren$, $N\ge 1$, with exponent $m>0$. The classical heat equation corresponds to the case $m=1$, while the range $m>1$ is called slow diffusion and the range $0<m<1$ is called fast diffusion. Extensive work has been done on \eqref{PME}, we mention as a first reference to  for $m\ne 1$ the book \cite{Vlibro} for the slow range, and \cite{VazSmooth} for the fast range. We will consider nonnegative solutions.

The combination of nonlinear diffusion of this type and strong anisotropy was recently considered by Feo, Volzone and the author  in \cite{FVV23}
by studying the \emph{anisotropic porous medium equation}
\begin{equation}\label{APM}
u_t=\sum_{i=1}^N(u^{m_i})_{x_i x_i}\quad \mbox{in  } \ \quad Q:=\mathbb{R}^N\times(0,+\infty)
\end{equation}
with $N\geq2$, and $m_i>0$ for $i=1,...,N$. The study was done in the fast diffusion range $0<m_i<1$ for all $i=1,\dots,N$.
We considered solutions to the Cauchy problem to \eqref{APM} with nonnegative and integrable initial data
\bigskip
\begin{equation}\label{IC} u(x,0)=u_0(x), \quad x\in \R^{N}.
\end{equation}
and  $u_0\in L^1(\mathbb{R}^N)$, $u_0\ge 0$. We put $M:=\int_{\mathbb{R}^N} u_0(x)\,dx$, so-called total mass. As in \cite{FVV23} we will focus here on the fast diffusion range that has some favorable features.  We need some assumptions on the exponents. We set
\begin{equation}\label{m bar}
m_c:=1-\frac{2}{N}, \quad \overline{m}:=\frac1{N}{\sum_{i=1}m_i}.
\end{equation}
We recall that  $m_{c}$ is the \emph{critical exponent} for the isotropic fast diffusion equation (\emph{i.e.} equation \eqref{APM} with $m_{1}=m_{2}=...=m_{N}=m<1$).
In this first work we will always assume the conditions

\noindent  \qquad (H1) \qquad $0< m_i< 1$ \qquad  for all $i=1,...,N$,

\noindent \qquad  (H2) \qquad $\overline{m}>m_c,$ i.e., $\sum_{i=1}m_i>N-2$.

Condition (H2) was essential for our results, being a necessary and sufficient condition in the isotropic case for the existence of fundamental solutions.  On the other hand, (H1)  is a condition of {\sl fast diffusion in all directions} that is made here for convenience, it is necessary for the analysis that we do. The case where there are slow diffusion exponents $m_i>1$ is not considered here. What we do below in the present paper does not apply to it.

Note that both basic conditions imply that $m_i<\overline{m} +2/N$, a bound that has an interest in the control of the amplitude of the anisotropy, a quality that is usually measured in terms of the differences $m_i-m_j$.

\medskip

\noindent {\bf Brief outline. } Paper \cite{FVV23} contains the construction of a theory of existence and uniqueness of continuous weak solutions, the properties of the ensuing evolution semigroup (for instance it is a semigroup of $L^1$-contractions with bounded orbits); plus the study of the existence and uniqueness of a self-similar fundamental solution, FSS for short, for every mass value along with its main properties, like positivity, boundedness and decay; and finally the asymptotic behaviour as $t\to \infty$ of all nonnegative solutions to the Cauchy Problem with the above conditions on data and exponents.

In this paper we want to extend the existing theory by introducing a  special singular solution of the type called the very singular solution, VSS for short. The VSS is a solution that, starting from an infinite mass located at the origin of coordinates as
initial datum, evolves according to the AFDE as an admissible solution away from the fixed singularity. The VSS shows the effect of anisotropy in a very quantitative way. In fact, we will show that the anisotropic diffusion distributes mass from the initial infinite mass singularity  along the different space directions: up to constant factors, there is a simple partition formula for this anisotropic mass expansion, as the minimum of separate 1-D VSS solutions. This is as a consequence of the strong scaling properties of the VSS and it has immediate consequences for many other solutions.

Once this object is analyzed, we will be able to control the exact behavior of the FS's both as $|x|\to  \infty$ and $t\to\infty$.
Actually, their tail is asymptotically convergent to the unique VSS tail, so there is a VSS partition formula for their diffusion at spatial infinity.

With this information we improve the tail analysis of general solutions with nonnegative initial data with finite mass. We obtain as consequences two important qualitative results: the Global Harnack Principle, GHP, and the asymptotic convergence in relative error, ACRE. The validity of these two properties for different nonlinear parabolic models is the object of much current investigation. We recall they do not hold in the same clean formulation for the heat equation or the porous medium equation posed in $\ren$.

We devote Section \ref{sec.hist} to review the relevant literature and some pertinent historical notes about VSS, GHP, asymptotic convergence, and the finer concept of ACRE. Section \ref{sec.op} concludes the paper with a commented list of extensions and open problems.

Lastly, our results  in dimensions $N\ge 2$ match perfectly with what is known in the case $N=1$ and in the isotropic version in $N\ge 2$ where these topics well studied. Our presentation offers a new approach and some new results for those topics.


\subsection{Self-similar solutions}\label{sec.sss}

We present next the main facts of the self-similarity analysis as done in \cite{FVV23}, that will be used here. Self-similar solutions of equation \eqref{APM} have the form
\begin{equation}\label{sss}
U(x,t)=t^{-\alpha}F(t^{-a_1}x_1,..,t^{-a_N}x_N)
\end{equation}
with constants $\alpha>0$, $a_1,..,a_n\ge 0$ to be chosen below. We substitute this formula into equation \eqref{APM}.
\noindent Note that  writing $y=(y_i)$ and $y_i=x_i \,t^{-a_i}$, we have
$$
U_t=-t^{-\alpha-1}\left[\alpha F(y)+
\sum_{i=1}^{N}a_iy_i\,F_{y_i}
\right],
$$
and
$$
\sum_{i=1}^{N}(U^{m_i})_{x_ix_i}=\sum_{i=1}^{N}t^{-(\alpha m_i+2a_i)}(F^{m_i})_{y_iy_i},
$$
Therefore, equation \eqref{APM} becomes
$$-t^{-\alpha-1}\left[\alpha F(y)+
\sum_{i=1}^{N} a_iy_i\,F_{y_i}
\right]=\sum_{i=1}^{N}t^{-(\alpha m_i+2a_i)}(F^{m_i})_{y_iy_i}.
$$
We see that time is eliminated as a factor in the resulting equation on the condition that:
\begin{equation*}\label{ab}
\alpha(m_i-1)+2a_i=1 \quad \mbox{ for all } i=1,2,\cdots,  N.
\end{equation*}
We also want integrable solutions that will enjoy the mass conservation property, which
implies $\alpha=\sum_{i=1}^{N}a_i$.
Imposing both conditions, and putting $a_i=\sigma_i \alpha,$ we get unique values for $\alpha$ and $\sigma_i$:
\begin{equation}\label{alfa}
\alpha=\frac{N}{N(\overline{m}-1)+2},
\end{equation}
and
\begin{equation}\label{ai} \ 
 \sigma_i= \frac{1}{N}+ \frac{\overline{m}-m_i}{2},
\end{equation}
so that \ $\sum_{i=1}^{N}\sigma_i=1. $ Observe that by Condition (H2) we have $\alpha>0$,  so that the self-similar solution will decay in time in maximum value like a power of time. This is a typical feature of a diffusion process.

As for the $\sigma_i $ exponents, we know that \ $\sum_{i=1}^{N}\sigma_i=1$, and in particular $\sigma_1=1/N$ in the isotropic case.  Condition (H1)  on the $m_i$  ensures that $\sigma_i> 0$. This means that the self-similar solution expands  as time passes, or at least does not contract, along any of the space variables.

With these choices,  the profile function $F(y)$ must satisfy the following nonlinear anisotropic stationary equation in $\mathbb{R}^N$:
\begin{equation}\label{StatEq}\sum_{i=1}^{N}\left[(F^{m_i})_{y_iy_i}+\alpha \sigma_i\left( y_i F\right)_{y_i}
\right]=0.\end{equation}
Conservation of mass must also hold : $\int U(x,t)\, dx =\int F(y)\, dy=M<\infty $ \ {\rm for} $t>0$.

 It is our purpose to prove that there exists a suitable solution of this elliptic equation, which is the anisotropic version of the equation of the Barenblatt profiles in the standard PME/FDE, cf. \cite{Vascppme, Vlibro}. \nlc
 The solution is indeed explicit in the isotropic case:
 $$
 F(y;m)=\left( C + \frac{\alpha(1-m)}{2mN}|y|^2\right)^{-2/(1-m)},
 $$
 with a free constant $C>0$ that fixes the total mass of the solution, $C=C(M)$. It is clear that this formula breaks down for $m\le m_c$ (called very fast diffusion range), where many new developments occur, see the monographs \cite{VazSmooth} and \cite{BBDGV}. We do not get any explicit formula for $F$ in the anisotropic case, we only have suitable estimates, in particular regularity and decay. Anisotropy will be evident in the graphics of the level lines.

\subsection{Self-similar variables}
In several instances in the sequel it will be convenient to pass to self-similar variables,
by zooming the original solution according to the self-similar exponents \eqref{alfa}-\eqref{ai}. More
precisely, the change is done by the formulas
 \begin{equation}\label{NewVariables}
v(y,\tau)=(t+t_0)^\alpha u(x,t),\quad \tau=\log (t+t_0),\quad y_i=x_i(t+t_0)^{-\sigma_i\alpha} \quad i=1,..,N
\end{equation}
with $\alpha$ and $\sigma_i$ as before.
We recall that all of them are positive. There is a free time parameter $t_0\ge 0$ (a time shift).

\begin{lemma}\label{Lem1}
If $u(x,t)$ is a solution (resp. super-solution, sub-solution) of \eqref{APM}, then $v(y,\tau)$ is a solution (resp. super-solution, sub-solution) of
\begin{equation}\label{APMs}
v_\tau=\sum_{i=1}^N\left[(v^{m_i})_{y_i y_i}+\alpha \sigma_i \left(\,y_i\,v\right)_{y_i}\right] \quad \mathbb{R}^N\times(\tau_0,+\infty).
\end{equation}
\end{lemma}
This equation will be a key tool in our study. Note that the rescaled equation does not change with the time-shift  $t_0$ but the initial value in the new time does, $\tau_0 = \log(t_0)$. If $t_0 = 0$ then $\tau_0 =-\infty$ and the $v$ equation is defined for $\tau\in \mathbb{R}$.

We stress that this change of variables preserves the $L^1$ norm. The mass of the $v$ solution at new time $\tau\geq\tau_0$ equals that of the $u$ at the corresponding time $t\geq0$.


\section{Preliminaries on the basic theory}\label{sec.bt}

The basic existence and uniqueness of solutions of the Cauchy Problem for equation  \eqref{APM} is stated in  Theorem 2.1  of \cite{FVV23} in that paper in the class of non-negative and integrable initial data, and proved in detail along with many useful properties.

 \begin{theorem}\label{EUWES}  Let the exponents $m_i$ satisfy assumptions (H1) and (H2). Then, for any nonnegative $u_0 \in L^1(\mathbb{R}^N)$ {there is a unique function $u\in C([0,\infty): L^1(\mathbb{R}^N))$ such that $u, u^{m_{i}}\in L^{1}_{loc}(Q)$ for all $i=1,...,N$, and equation \eqref{APM} holds in the distributional  sense in $Q=\mathbb{R}^N\times(0,+\infty)$, with the following additional properties:}

1)  $u(x,t) $ is a uniformly bounded function for each $\tau>0$ with an estimate of the form $\|u(\cdot,t)\|_\infty\leq C t^{-\alpha}$. The precise estimate is given by
\begin{equation}\label{smooth}
\|u(t)\|_\infty\leq C t^{-\alpha}\|u_0\|_1^{2\alpha/N}\quad \forall t>0.
\end{equation}

2)  Let $Q_\tau=\ren\times (\tau,\infty)$. We have $\partial_i u^{m_{i}}\in L^{2}(Q_\tau)$ for every $i$ and the energy estimates are satisfied:
\begin{equation*}\label{Energywholespace}
\int_0^T\int_{\ren}
\left|\frac{\partial}{\partial x_i}u^{m_i}\right|^2 \, dx\, dt
\leq\int_{\ren}\left[
\frac{1}{m_i+1}{u_0}^{m_i+1}
\right]\, dx
-
\int_{\ren}\left[
\frac{1}{m_i+1}{u}^{m_i+1}(x,T)
\right]\, dx
\end{equation*}
for all $i=1,...,N$ and $T>0$.  Equation \eqref{APM} holds in the usual weak sense applied in $Q_\tau$ for every $\tau>0$.

3) Consequently, the maps $S_t: u_0\mapsto u(\cdot,t)$ generate a semigroup of $L^1$ ordered contractions in $L^1_+(\ren)$. The $L^1$-contraction estimates are satisfied. The maximum principle applies.

4) Conservation of mass holds: for all $t>0$ we have \ $\int u(x,t)\,dx=\int u_0(x)\,dx$. Assumption (H2) is crucial for mass conservation.

5) If we start with initial data $u_0\in L^1(\ren)\cap L^\infty(\ren)$ we may also conclude item 2) with $\tau=0$ and $u(x,t)$ is uniformly bounded and continuous  in space and time.
 \end{theorem}

A further topic that is settled in that paper  is the existence and properties of finite-mass self-similar solutions,  solved in Theorem 1.1 of  \cite{FVV23}.

\begin{theorem}\label{fundamental solution} Under the restrictions (H1) and (H2), for any mass $M>0$ there is a unique self-similar fundamental solution $U_M(x,t)\ge 0$ of equation \eqref{APM} with mass $M$. The profile $F_M$  of such a solution is an SSNI (separately symmetric and nonincreasing) positive function. Moreover,   $F_M(y)$ is  positive everywhere, bounded and smooth, and decays at infinity in a power-like way.
\end{theorem}

Note that uniqueness is proved in the class of self-similar mass conserving solutions. A further study of such solutions will be done below in the present paper by relating them to the study of more singular object called the Very Singular Solution, whose simplicity is due to its increased invariance properties.

 A further result of \cite{FVV23} establishes the asymptotic behaviour of finite mass solutions in Theorem \ref{EUWES}.

\begin{theorem}\label{thmasympto}
Let $u(x,t)$ be the unique solution  of the Cauchy problem for equation \eqref{APM} with nonnegative initial data $u_{0}\in L^{1}(\R^{N})$  under the restrictions (H1) and (H2). Let $U_{M}$ be the unique self-similar fundamental \nlc solution with the same mass as $u_{0}$. Then,
\begin{equation}\label{L1conv}
\lim_{t\rightarrow\infty}\|u(t)-U_{M}(t)\|_{L^1(\ren)}=0.
\end{equation}
The convergence holds in the $L^{p}$ norms, $1\le p <\infty$,  in the proper scale
\begin{equation}\label{Lpconv}
\lim_{t\rightarrow\infty}t^{\frac{(p-1)\alpha}{p}}\|u(t)-U_{M}(t)\|_{L^p(\ren)}=0,
\end{equation}
where $\alpha=N/(N(\overline{m}-1)+2)$ is the constant in \eqref{alfa}.
\end{theorem}

Under certain conditions of the initial data, convergence \eqref{Lpconv} is proved also for $p=\infty$. But the results in this case are not very sharp and we will devote part of our effort to prove a improved  versions in the the form of GHP and relative error convergence, see Theorems \ref{GHP} and \ref{ACRE}.


\section{The local mass estimate }\label{local.mass}

It can  be called the Anisotropic Herrero-Pierre estimate, since the isotropic version of the estimate was  was introduced by them for the FDE in \cite{HP}. It works in the isotropic on the condition that  $0<m_i<1$ for all $i$. Then  there is a local estimate for nonnegative solutions of the AFDE as in the isotropic FDE.

We start from the equation \eqref{APM}, multiply it by a smooth test function $\vp(x)\ge 0$ with compact support, and integrate to get
$$
\frac{d}{dt}\int_{\ren} u\,\vp\,dx=\sum_i \int_{\ren}  u^{m_i}\partial^2_{ii}\vp\,dx=\sum_i \int_{\ren}  u^{m_i}\vp^{m_i}\frac{\partial^2_{ii}\vp}{\vp^{m_i}}\,dx
$$
Let us call $X(t)=\int_{\ren}  u\,\vp\,dx$. Then, using H\"older we get
$$
\frac{dX}{dt}\le \sum_i C_i X(t)^{m_i}\, Y_i^{1-m_i}, \quad
Y_i=\int_{\ren}  \left(\vp^{-m_i}\partial^2_{ii}\vp\right)^{1/(1-m_i)}\,dx
$$
Note that the $Y_i$ are numbers that depend only on $\vp$. We thus get an ODE inequality that allows the mass to be controlled locally as in the isotropic FDE.

\medskip

We want to estimate $Y_i$ in a useful way by choosing a suitable $\vp$. If we assume that $\vp $ is supported in a compact set $K\subset \ren$, we may improve the formula for $Y_i$ into
$$
Y_i=\int_K \left(\vp^{-m_i}\partial^2_{ii}\vp\right)^{1/(1-m_i)}\,dx_1\cdots dx_N.
$$
Now we take $\vp $ a product of separate functions:  $\vp(x)=\Pi_i \vp_i(x_i)$,  and assume that the value of each $\vp_i(x_i)$  approaches the border of the supporting interval $K_i=[a_i, b_i]$ in a very flat way that we will determine later. The integrals in the definition of the $Y_i$ must be finite. We normalized in lengths by putting $\vp_i(x_i) = \phi_i(t)$ with $x_i = a_i + L_i t$, $0 < t < 1$. Let us then calculate the first integral as an example of what happens. Since
$$
\partial^2_{11}\vp(x)= \partial^2_{11}\vp_1(x_1) \,\Pi_{j\ne 1}\vp_j(x_j),
$$
inserting it into the formula for $Y_1$ and integrating separately we get
$$
Y_1= c_1 Z_1 L_2\dots L_N, \qquad
Z_1 =\int_{K_1} \left(\vp_1^{-m_1}\partial^2_{11}\vp_1\right)^{1/(1-m_i)}\,dx_1.
$$
The last integral is well known to be finite if the $\vp_1$  approaches $b_1$ like $C(b_1-x_1)^k$, $k>1$ is large, and same rate at $x_1$ approaches $a_1$. Indeed, $\partial_{1}\vp_i(x_i)$ decays like   $C(b_1-x_1)^{k-1}$. Then, $\partial^2_{11}\vp_1$ is bounded and
$$
\partial^2_{11}\vp_1 \sim C(b_1-x_1)^{k-2}
$$
Hence, $\vp_1^{-m_1}\partial^2_{11}\vp\sim C(b_1-x_1)^{k-2-m_1 k}$, so that the last integral of $Y_1$ is a bounded number if $k(1-m_1)>1$.
If we want to see the dependence of $Z_1$ on the length of the interval we put $\vp_1(x_i)=\phi(t)$ with $x_i=a_1+ L_1 t, $ $0<t<1$, and we have
$$
Z_1 =\int_{K_1} \left(\phi^{-m_1}(t)L_1^{-2}\partial^2_{11}\phi(t)\right)^{1/(1-m_i)}\,L_1dt=L_i^{1-(2/(1- m_1))}
$$
Finally, since the volume of $K$ is $V(K)=L_1\cdots L_N$, we get
\begin{equation}\label{alme.ode}
\left|\frac{dX(t)}{dt}\right|\le C_i \sum_i X(t)^{m_i} Y_i^{1-m_i}, \quad
Y_i=V(K)\,L_i^{-2/(1- m_i)} C_i.
\end{equation}
We get an ODE that allows the mass to be controlled locally as in the isotropic FDE.

\medskip

{\bf Renormalization. } \noindent
We may still may simplify the result by defining the normalized local mass as $\widetilde X(t)= X(t)/V(K)$.
then diving by $V(K)$ and putting
$$
\left|\frac{d \widetilde X(t)}{dt}\right|\le \sum_i C_i \widetilde X(t)^{m_i} \, L_i^{-2}.
$$

The solution of the ODE is not explicit. But, of course for large values of $ \widetilde X(t)$ the ODE will like
$$
Z'(t)\simeq c C_1 z(t)^{m_i}
$$
that the finite solution will look like $z^{1-m_1}(t)=z^{i-m_1}(0)+ D_1\,(t /L_1^{2})$ where $m_1$ is taken as the largest of the  $m_i$. This will be a supersolution,)  so that for all large $\widetilde X(t)$ we get
\begin{equation}\label{alme.upper}
\widetilde X(t)=O(t^{1/(1-m_1)}).
\end{equation}
In doing that simplification we will miss some detailed information, but it is enough for our purposes for the moment.


\medskip

\section{The  Very Singular Solution}\label{sec.vss}

We concentrate now on the construction of this object that will play a great role in the remaining theory.

\medskip

\noindent {\bf 1. Existence of the VSS as limit of fundamental solutions}

 The local boundedness estimate derived from \eqref{alme.ode} can be used to prove that  the self-similar increasing family of fundamental solutions $U_M$ admits a finite limit of when the mass $M\to \infty$:
$$
V_\infty(x,t)=\lim_{M\to \infty}U_M(x,t)<\infty.
$$
This limit  taken in a monotone increasing way and it is bounded  for $x\ne 0$ and all $t>0$. We call it the Very Singular Solution. It is still a self-similar solution of the same type as the $U_M$ with profile space $F_\infty(y)$
$$
V_\infty(x,t)=t^{-\alpha}F_\infty(...,x_i\,t^{-\sigma_i\alpha},...)
$$
with $\alpha$ and $\sigma_i$ as introduced in Section \ref{sec.intro}.

We can use the local ODE \eqref{alme.ode} to conclude that $F_\infty(0)=+\infty$ and $\int_B F_\infty(x)\,dx=+\infty$ in any ball near $x=0$. On the other hand, the tail of $ F_\infty$ has mass zero at infinite.

\medskip

\noindent {\bf  2. Anisotropic behaviour of the VSS}

The VSS inherits the scale invariance of all the $U_M. $ A key property of the solution $V_\infty$ is the new invariance under the change-of-mass rescaling. Indeed, if we recall that for every $k>0$ we have
\begin{equation}\label{mass.change}
F_{k^{\beta }M}(y_1,\cdots,y_N)=k F_M(k^{\gamma_i}y_i), \quad  , \ \gamma_i=(1-m_i)/2>0
\end{equation}
and some $\beta=1-N(1-\overline {m})/2>0 $. Letting $M\to\infty$ we get the new identity
$$
F_\infty(y_1,\cdots,y_N)=k F_\infty(k^{\gamma_i}y_i), \quad \gamma_i=(1-m_i)/2.
$$
This is crucial in what follows.

\medskip

\noindent {\bf  Behaviour on the unit sphere $ S=\mathbb{S}^{N-1}$.} We take $y=\omega$ of unit length and we get
a function $C(\omega)=F_\infty(\omega)$ that is bounded above and below away from 0, thus a $C^\infty$ function by the theory of local regularity
of locally non-degenerate solutions of parabolic equations.

\medskip

\noindent {\bf  Behaviour along the axes.} Now we take $y=r\,e_i$, $e_i=(...,1,...)$ being the unit vector in the $i$ direction. If we put $r=k^{-\gamma_i}$ in the above formula then for all $r>0$ we get  the explicit form
$$
F_\infty(r\,e_i)=C_i\, r^{-2/(1-m_i)}, \quad C_i=F_\infty(e_i).
$$
Going back to $V_\infty$ we get along every axis the explicit form
\begin{equation}
V_\infty(r\,e_i, t)=C_i\,t^{\frac{1}{1-m_i}}r^{-2/(1-m_i)},
\end{equation}
Note that the time exponent of the self-similar solution is given by
$$
 \delta_i=\frac2{1-m_i}\sigma_i\alpha-\alpha=\alpha( \frac{2\sigma_i}{1-m_i}-1)=\frac{1}{1-m_i}.
$$
It is the same form as in the isotropic case. It is a kind of ``separate behaviour" along the axes with a weak interaction in between the directions that we will describe next. The constant $C_i$ may depend on all the $m_i$'s. The strong likelihood looks like a miracle.

\medskip

\noindent {\bf  Formula for $F_\infty(y)$ in the general case.} For points that do not lie on the axes, we proceed as follows. We can completely cover $\ren\setminus \{0\}$ with lines of the form
$$
y_i(k,\omega)=k^{\gamma_i} \omega_i, \qquad \omega\in \mathbb{S}^{N-1}, \   k>0.
$$
If we  follow the lines $y(k,\omega)$ for $k>0$ with $\omega$ fixed, then
$$
F_\infty(y(k,\omega))=\frac{F_\infty(\omega)}{k}.
$$
We write this in the form

 \begin{theorem}\label{AVSS} for every $y_i\ne 0$ we have
\begin{equation}\label{vss.formula1}
F_\infty(y(k,\omega))= C(\omega) \,|\omega_i|^{2/(1-m_i)}\,|y_i|^{-2/(1-m_i)},
\end{equation}
 where $C(\omega)=F_\infty(\omega)>0$ is a smooth function. As a consequence we get the following anisotropic inequalities: There are positive constants constants $K_1<K_2$ (depending on $N$ and the exponents $m_i$) such that
\begin{equation}\label{vss.formula1b}
\frac{K_1}{  \sum_i |y_i|^{2\mu_i}} \le
F_\infty(y)\le \frac{K_2}{  \sum_i |y_i|^{2\mu_i}}, \qquad \mu_i=\frac1{1-m_i}>1.
\end{equation}
\end{theorem}

Formula \eqref{vss.formula1} completely determines $F$ along the covering of the whole space by lines that go from $y=0$ to infinity allows in an easy way to describe the \sl level lines \rm $F={L}>0$ in the parametric form
$$
x_i (\omega, L)=\omega_i \,C(\omega)^{\gamma_i}{L}^{-\gamma_i}.
$$
Let us now get the inequalities \eqref{vss.formula1b} and obtain a better understanding of the decay estimate.

\noindent {\sc  Bound from above:} We see that for every $y\ne 0 $ there exists a coordinate $y_i$ that is not zero, assume it is positive. Then we have
$$
F_\infty(y)=C(\omega_i)\,(|\omega_i|/|y_i|)^{-2\mu_i}\le C_{M}|y_i|^{-2\mu_i},
$$
where $\mu_i=1/(1-m_i)$ and $C_M=\max\{C(\omega), |\omega|=1 \}$. Since the bound holds for all $i$ such that $y_i\ne 0$ we get
\begin{equation}\label{vss.formula2}
F_\infty(y)\le C_{M}\min_i\{ |y_i|^{-2\mu_i}\}.
\end{equation}
Besides, using the arithmetic inequality for positive numbers
$$
\max(a_1,\cdots,a_N)\ge \frac1{N} (a_1+...a_N),
$$
and applying it to $a_i= |y_i|^{2\mu_i}$, we arrive at the form
\begin{equation*}\label{vss.formula3}
F_\infty(y)\le N C_{M}\left( \sum_i |y_i|^{2\mu_i}\right)^{-1}.
\end{equation*}
These are the desired clean forms for the upper bound.

\medskip

\noindent {\sc  Bound from below:} To see that the obtained bound is sharp, we see that same argument works when we replace $C_{M}$ by $C_{min}$ and we get
$$
F_\infty(y)\ge C_{min} \min_i \{ \omega_i^{2\mu_i}r_i^{-2\mu_i}\}
$$
Given $y$ with positive coordinates, let for instance $i=1$ the exponent that realizes the maximum of $\{|y_i|^{\mu_i}: i=1,\cdot, N\}$. Then we have for every $i\ne 1$ along the curve that joins $y$ to $\mathbb{S}^{N-1}$
$$
\omega_i^{2\mu_i}/\omega_1^{2\mu_1}=|y_i|^{2\mu_i}/|y_i|^{2\mu_1}\le 1
$$
which means that $\omega_1^{2\mu_1}$ attains the maximum value among the $\omega_i^{2\mu_i}$. Since $\sum_i \omega_i^2=1$ it follows has it has to be larger than some $c_0>0$ that depends on the exponents at hand.
In that case
$$
F_\infty(y)=C(\omega)/k=C(\omega) |\omega_1|^{2\mu_1}|y_1|^{-2\mu_1} \ge C_* |y_1|^{-2\mu_1}\ge C_{*} \min_i \{ |y_i|^{-2\mu_i}\}
$$
and since we have the arithmetic inequality:  $\max(a,b)\le a+b$, we get
\begin{equation*}\label{vss.formula4}
F_\infty(y) \ge C_{*} \left(\sum_i |y_i|^{2\mu_i}\right)^{-1}.
\end{equation*}

\medskip

\noindent   {\bf 3. Comparison with separate 1D diffusions }

If we make the separate analysis in the 1-D case with some exponent $m\in (0,1)$ we get the fundamental solutions of the form
\begin{equation}\label{vss.formula5}
 F_{M,1}(y;m)=\left( C + \frac{(1-m)}{2(1+m)}\,|y|^{{2}}\right)^{{{-1}}/(1-m)},
\end{equation}
 with a free constant $C>0$ that fixes the total mass  $M$ of the solution, $C=C(M)$. Letting $C\to 0$ we get the
 1-D VSS with profile
 $$
  F_{\infty,1}(y;m)= C(m; 1)\,|y|^{-{2}/(1-m)},\qquad C(m;1)= \left(\frac{(1-m)}{2(1+m)}\right)^{{{-1}}/(1-m)}\,.
 $$
 Putting $m=m_i$ this is exactly the expression of the terms of the $N$-dimensional formulas \eqref{vss.formula1}-\eqref{vss.formula1b} but for a multiplicative constant  that depends on the set of $m_i$'s. The complete time dependent 1-D VSS is
\begin{equation}\label{vss.1d}
 V_{\infty,1}(x,t;m)= C(m; 1)\,t^{1/(1-m)} |x|^{-{2}/(1-m)}.
 \end{equation}
 We may interpret formulas \eqref{vss.formula1} to \eqref{vss.formula1b} as the way the matter that  spreads around in anisotropic diffusion is distributed  along the different directions. Indeed, all directions share in the resulting solution according to the predicted 1D strength of each of the directions, and the joint formula is a fair partition  up to controllable constants. Later we will examine how much of this universal fair partition is conserved when we deal with fundamental solutions or more general classes of solutions.

\medskip

\noindent {\bf  Example in 2D.} Note that $m_c=(N-2)/N=0$, We that $1>m_1>m_2>0$, for instance $m_1$ close to 1 and $m_2$ near zero.
The behaviour  along the axes is
$$
V_\infty(r_1,0)=C_1\,(t/r^2_1)^{\mu_1}, \qquad V_\infty(0,r^2_2)=C_2\,(t/r_2^2)^{\mu_2},
$$
with $\mu_i= \frac{1}{1-m_i}$. Therefore, $\mu_1\approx \infty$ which means  very fast decay along the first axis, and $\mu_2\approx 1, $ close to the smallest possible, which means a very fat tail in that direction. $F_\infty$ will be given by the value $C/k$ along lines where $k$ varies for fixed $\omega=(\omega_1, \omega_2)$ which is given by
$$
y_1=\omega_1\, k^{\gamma_1}, \quad y_2=\omega_2\, k^{\gamma_2},
$$
where $\gamma_i= 1/(2\mu_i)$, which gives the line $y_2= K(\omega)y_1^{\mu_1/\mu_2}$. It looks like a parabola looking at the second axis. We conclude that the fastest decay is dominant in the space as we approach infinite.

 Note that whenever $1> m_1>m_2$ then
$$
\sigma_2-\sigma_1= \frac12 (m_1-m_2),
$$
hence the second axis expands more than the first, the $y(k)$ curve looks to the second axis. Moreover,
$$
\gamma_1=\frac1{1-m_1}>\gamma_2=\frac1{1-m_2}
$$
so that the first axis decays more than along the second. The faster decay is dominant in the space.

\medskip

\noindent   {\bf 4. Monotonicity in time of the VSS and time stability}

 We look at the multidimensional anisotropic case. We want to show that $V_t\ge 0$ as in the 1D case. In view of the self-similar we have
  $$
\partial_t V(x,t)= -\alpha \,t^{-(\alpha+1)} \left(F(y) +\sum_i \sigma_i y_i F_{y_i}\right),
$$
hence we need to prove that
$$
F(y) +\sum_i \sigma_i y_i F_{y_i}\le 0 \quad \forall y\ne 0.
 $$
 This is true on the axes where only one $y_i$ is nonzero and there $F(y)=C_i\,|y_i|^{-2/(1-m_i)}$, so that
 $$
 F(y) +\sum_j \sigma_j y_j F_{y_j}=F -\frac{2\sigma_i}{1-m_i}\,F =-\frac{1}{\alpha(1-m_i)}\,F<0\,.
 $$
In the general situation, we move along the line $y(k,e)$ and then we have $F=C(\omega)/k$, hence
$$
\frac{\partial F}{\partial k}=-\frac1{k^2}C(\omega)=-\frac{F}{k}, \quad \frac{\partial F}{\partial k}= \sum_i F_{y_i}\frac{\partial y_i}{\partial k}= \sum_i \gamma_i\frac{y_i}{k}F_{y_i}
$$
so that $ F=-\sum_i \gamma_i y_i\,F_{y_i}. $
Hence,
$$
\partial_t V = \alpha \,t^{-(\alpha+1)}\,(\sum_i \gamma_i y_i\,F_{y_i} -\sum_i\sigma_i y_i F_{y_i})=
 \alpha \,t^{-(\alpha+1)}\, \sum_i(\sigma_i- \gamma_i)|y_i| |F_{y_i}|.
$$
Now $\sigma_i- \gamma_i=1/(2\alpha)$, so that
$$
\partial_t V =  \frac12 \,t^{-(\alpha+1)}\, \sum_i |y_i| |F_{y_i}|.
$$
this proves that $V$ is increasing in time. But we have a finer estimate in time. If $1>m_1>m_2>...>m_N$ we get the reverse order for the gammas and then
$$
\sum_i |y_i| |F_{y_i}|\le \frac1{\gamma_1} \sum_i \gamma_i|y_i| |F_{y_i}=\frac1{\gamma_N}F
$$
so that in the end
$$
\mu_N \le  \frac{t\partial_t V}{V} \le \mu_1.
$$
These estimates are sharp. In particular we can calculate the relative error of a delay in time as follows:

\begin{corollary} \sl For all  $t,h>0$  we have
\begin{equation}\label{asymp.V}
\frac{V_\infty(x,t+h)-V_\infty(x,t)}{h\,V_\infty(x,t)}\approx C \frac{\partial_t V_\infty(x,t)}{V_\infty(x,t)}\approx O(\frac{1}{t}).
\end{equation}
The approximate signs apply for $h$ much smaller than $t$.
Thus, the relative error quotient goes to zero uniformly for $t\to\infty$ with rate $O({t}^{-1})$.
\end{corollary}

We immediately see that in the case of time delayed solution $u(x,t)=V_\infty(x,t\pm h)$, the relative convergence rate $O({t}^{-1})$ is optimal. Extending this asymptotic estimate to more general solutions is an important task. Note that similar rates appear in the isotropic theory in \cite{CV03} and ensuing literature.

\medskip

\noindent {\bf Space translations.} We may ask if the estimate of the relative error applies to the singular solutions obtained by space translation, $u(x,t)==V_\infty(x-a,t\pm h)$. It can be shown that estimates in relative error happens as $t\to\infty$ only on outer sets that recede away from the origin.


\section{Universal tails and  Global Harnack Principle}\label{sec.GHP}

The precise knowledge of the VSS allows to derive interesting consequences for the fundamental solutions $U_M$. Thus,
immediate comparison shows that the FS $U_M$ has a good a priori upper bound everywhere
$$
U_M(x,t)\le V_\infty(x,t)\le C\min_i\{(t/|x_i|^2)^{\mu_i}\}     \qquad \forall x, \ \forall t>0.
$$
A corresponding lower bound and tail behaviour for the fundamental solutions is done by rescaling the fundamental profile. Since $F_M$ that tends to $F_\infty$ locally uniformly as  $M\to\infty$, rescaling turns  into the equivalent
statement
$$
\lim_{|y|\to\infty}\frac{|F_M(y)-F_\infty(y)|}{F_\infty(y)}=1,
$$
where we use formula \eqref{mass.change}. This is the universal tail behaviour we were aiming at: all self-similar solutions exactly copy the tail behaviour  from $F_\infty$ in the very precise relative error. Equivalently, we get for every $t>0$
$$
\lim_{|x|\to\infty}\frac{|U_M(x,t)-V_\infty(x,t)|}{V_\infty(x,t)}=1,
$$
and this limit is uniform in outer sets of the form $A(C)=\{(x,t): |x_i|\ge C\, t^{\sigma_i \alpha}\}$.

The universal tail behaviour of the fundamental solutions that we have proved is  a key to obtain the
 result about Global Harnack Principle.

\begin{theorem}\label{GHP} Let $u$ be a solution of our Cauchy Problem such that the tail of the initial datum lies below the space profile of a VSS  at some time $t_1$. Then for every time $T>0$ there are constants $0<C_1<C_2< \infty$ such that
\begin{equation}\label{ghp.1}
C_1 U_M(x,t)\le u(x,t)\le C_2 U_M(x,t)
\end{equation}
for all $x\in\ren$ and $t\ge T$. The constants $C_1<1<C_2$ depend on $T$ and the data $u_0$.
\end{theorem}

The proof of the theorem is divided into upper and lower estimates.

\noindent  {\sc Lower bound.} We need some control on the nontrivial data. So we first recall that for every $t_1>0$ the solution is positive in any ball, see \cite{FVV23}, let us say that
$$
u(x,t)\ge c>0 \ \qquad \mbox{for } \  x\in B_R(0), \ t_1/2\le t\le t_1.
$$
Next, we make a comparison in the outer cylinder $Q=(\ren\setminus B_R)\times (0, t_1/2)$. More precisely, we compare $u(x,t+ (t_1/2))$ with the $U_{M'}(x,t)$ for some small mass $M'$. On the initial time line the comparison says that $u(x,t_1/2)\ge 0=U_{M'}(x,0)$. On the lateral boundary $u(x, t+(t_1/2))$ is bounded below away from zero by some $c>0$, while $U_{M'}(x,t)$ is smaller if $M'$ is small (with $t_1$ and $R$ fixed).
Hence,
$$
u(x,t_1)\ge U_{M'}(x,t_1/2) \ \qquad \mbox{for } \  x\in B_R(0).
$$
On the other hand, the same is true inside $B_R$ if $M'$ is small. We conclude from the MP that for all later times the two solutions are ordered:
$$
u(x,t+t_1)\ge U_{M'}(x,t+t_1/2) \ \qquad \mbox{for } \  x\in  \ren, \  t\ge 0.
$$
At this moment we only need to prove that  there exists $K_1 >0 $ such that
$$
U_{M'}(x,s)\ge K_1 U_{M'}(x,s+t_1/2) \qquad \mbox{for } \  x\in \ren, s\ge t_1/2.
$$
This relation is easy to check. Of course, $K_1$ depends on $M'$ and $t_1$ and the solution under consideration. This ends the proof of the lower bound.

\medskip

\noindent  {\sc Upper bound.}   Take  a time $T>0$. Here we begin by the inner behaviour of $u(x,T)$ where we can use the smoothing effect \eqref{smooth} that was proved in \cite{FVV23}. The upper comparison \eqref{ghp.1} is true there on an inner set that expands in space with the self-similar rates $t^{\sigma_i\alpha}$ in each direction. This means that  for any choice of $K>0$ it holds for  $\{(x,t) : |x_i|\le K \,t^{\sigma_i\alpha}, \forall i; t\ge T/2\}$.

In order to estimate the outer behaviour of the solution at $T$ we may use separate one-dimensional VSS barriers, $V_{\infty,1}(x_i,t;m_i)$  that we studied in Section \ref{sec.vss},   formula \eqref{vss.1d}. More precisely, if the estimate
$$
u_0(x)\le V_\infty(x,t_1)
$$
holds outside of the ball of radius $R_0>0$ we use as comparison function
$$
\widetilde V(x,0)=V_{\infty,1}(x-R_0,t_1)=C(m_1;1)\,t_1^{1/(1-m_1)}|x_1-R_0|^{-2/(1-m_1)}.
$$
We have $u_0(x)\le V_\infty(x-R_0,t_1)$ in the half-space $H_1=\{x\in \ren: \ x_1>R_0\}$. We may use the Maximum Principle on the two solutions defined in $Q_1=H_1\times [0,\infty)$ and conclude that
$$
u(x,t)\le \widetilde V(x,t)=V_{\infty,1}(x-R_0,t+t_1)= C\,(t+t_1)^{1/(1-m_1)}|x_1-R_0|^{-2/(1-m_1)},
$$
valid for all $x_1>R_0$ and $t>0$. Same formula for $x_1<-R_0$. A similar upper bound holds in the rest of the axis directions. If we compare this with formula \eqref{vss.formula4} and use the self-similarity to pass from $F_\infty$ to $V_\infty$ and consider only $t>T$ with $T$ large enough, we get the desired upper bound in the outer region not covered by the inner upper estimate.
\quad \qed  \nc


\section{Large-time behaviour of general solutions}\label{sec.AB}

Next, we give the result about fine asymptotics as $t\to \infty$. We prove ACRE, asymptotic convergence in relative error.

\begin{theorem}\label{ACRE} Let $u$ be a solution of our Cauchy Problem such that the tail of the initial datum lies below the profile of the VSS taken at some time $t_1>0$. Then we have
\begin{equation}\label{asym.acre}
\lim_{|t|\to\infty}\frac{|u(x,t)-U_M(x,t)|}{U_M(x,t)}=1.
\end{equation}
as a uniform limit in $x\in\ren$.
\end{theorem}

\noindent {\sc Proof.} (i) For the behaviour in the typical inner region, we accept the asymptotic theorem of \cite{FVV23} that ensures convergence of $u$ to the fundamental solution $U_M$  inside the expanding core of the domain \
$$
{\mathcal C}(k,T)= \{(x,t): t\ge T, \ |x_i|\le k\, t^{\sigma_i\alpha} \ \forall i \}
$$
in  relative sup norm. The word relative is justified by looking at the uniform convergence in fixed ball of the renormalized variable $v(y,\tau)$ defined in \eqref{NewVariables}. So, given $k>0$ and  $\ve$, we take an initial datum $u(x,T)$ which is $\ve$ near $U_M(x,T)$ in ${\mathcal C}(k)$ for  large $T\ge T(k,\ve)>0$. That shows the asymptotic result     \eqref{asym.acre} in that kind of inner region.

\medskip

(ii) For the outer behaviour we need bounds from above and below. The outer zone for $t\ge T$ will be denoted by
$$
{\mathcal O}(k,T)=\ren\times(T,\infty)\setminus {\mathcal C}(k)
$$
Now, from the proof of the GHP in previous section we know that  for $t_*>T$
$$
u(x,t_*)\ge U_{M'}(x,t_*- T/2),
$$
which follows from the outer behaviour of any fundamental solution. Hence, for very large $t_*$ we get in ${\mathcal O}(k,t_*$), if we restrict to $k$  large enough
$$
u(x,t)\ge (1-\ve)\,V_{\infty}(x,t- T/2)
$$
But for large $t_*$ and same outer region we then have
$$
u(x,t)\ge (1-2\ve)\,V_{\infty}(x,t)\ge (1-2\ve) U_{M}(x,t)
$$
with $t\ge t_*$, which is the desired lower bound in the outer zone.

\medskip

(iii) Finally, we look for an upper bound in the outer zone. We use again the fact that the initial tail is bounded above by a VSS.
If we look at the proof of the GHP in Theorem \eqref{GHP} we see that for all $|x_i|>R_0$ and $t>0$ we have
$$
u(x,t)\le C(m_i;1)\,(t+T_i)^{1/(1-m_i)}|x_i-R_0|^{-2/(1-m_i)},
$$
so that in the outer zone ${\mathcal O}(k)$ and for $t_*$ large enough
$$
u(x,t_*)\le \min_i \{ C(m_i;1)\,t_*^{1/(1-m_i)}|x_i|^{-2/(1-m_i)}\} \le C V_\infty(x,t_*)\le (1-\ve) V_\infty(t_*+ A).
$$
for some $A$. But in the outer zone we know that $(1-\ve) V_\infty(x,t_*+ A)\le  U_M(x,t_*+ A)$.
In this way we conclude that
$$
u(x,t_*)\le  U_M(x,t_*+ A) \qquad \forall  x  \in   \ren,
$$
which means by the maximum principle that for all $t>t_*$
$$
u(x,t)\le U_M(x,t+ A) \qquad \forall  x  \in   \ren.
$$
Finally, we use the outer behaviour of $U_M(x,t+ A)$. It is clear that
\begin{equation}\label{asym.acre}
\lim_{|t|\to\infty}\frac{|U_M(x,t+A)-U_M(x,t)|}{U_M(x,t)}=1.
\end{equation}
That ends the proof of the ACRE. \quad \qed
\nc

\section{Historical notes and references}\label{sec.hist}

\noindent    {$\bullet$} The reference for the previous theory of the anisotropic FDE \eqref{APM}   has been the paper with Feo and Volzone \cite{FVV23}. A similar $L^1$ theory  was written in parallel by the same authors in the work \cite{FVV21} for the  fast diffusion evolution equation with anisotropic nonlinear gradient diffusion of $p$-Laplacian type. Interesting previous references for the anisotropic FDE are \cite{H, SJ05, SJ06}.

\smallskip

\noindent    {$\bullet$}  The  very singular solution of the isotropic FDE posed in $\ren$ is given by the  formula
\begin{equation}\label{vss.iso}
V_\infty(x,t)=C(m,N)\,t^{1/(1-m)}|x|^{2/(1-m)},
\end{equation}
with an explicit constant $C(m,N)>0$. It is valid for $1>m>m_c$, and verifying the formula is immediate and it  was well known. Comparing with formulas \eqref{vss.formula1}--\eqref{vss.formula4} we may see that \eqref{vss.iso} is the particular form of our  anisotropic VSS solution when all the $m_i$ are equal. That isotropic  VSS was used in  \cite{ChV02} as minimum tail behaviour, a part of the study of solutions of the FDE with infinite data on some fixed sets.

Very singular solutions emanating from an infinite mass located at one point and taking a form like \eqref{vss.iso} have been constructed by us for the  a number of equations: the $p$-Laplacian equation, see \cite{VazSmooth}, for the fractional porous medium equation in the isotropic fast diffusion case \cite{Vaz2014}, and for the fractional $p$-Laplacian equation in \cite{Vaz20}. The similarity ideas of this last reference have strongly motivated the last part of the present paper.

Historically, there have been many works on the existence and role of VSS since the 1980's, mostly dealing with semilinear heat equations with reaction/absoption terms. Let us mention the influential paper by Brezis et al. \cite{BPT85} in 1985 that was followed by many interesting works.

\smallskip

\noindent    {$\bullet$}  The behaviour of the solutions of the FDE that leads to the Global Harnack Principle has been established in the isotropic FDE case in \cite{Vascppme}. The extension to solutions of the so-called weighted-FDE with sharp rates of convergence was done
in \cite{BonSim19} and \cite{BonSim23}, and extended to the isotropic fast diffusion $p$-Laplacian equation in  \cite{BonSimStan22}.
Similar work for the fractional $p$-Laplacian equation was done \cite{Vaz20}. The convergence in relative error ACRE is done in all those papers.
The present paper seems to be the first one to contain the two properties, GHP and ACRE, in the anisotropic setting.

\smallskip

\noindent    {$\bullet$} In the isotropic PME we know that no GHP or ACRE can be true in the stated data generality when $m\ge 1$, i.e., for linear diffusion or slow diffusion. However, it is known that GHP holds in that range in cases of corresponding fractional diffusion, though the existence of a VSS and the universal tail behaviour are not true.


\section{Comments,  extensions and open problems}\label{sec.op}

\noindent    {\bf Negative exponents.} Formally, the good exponent range: \ $\sum_{i=1}m_i>N-2$, $m_i<1$, allows us to consider exponents $m_i\le 0$ on the condition that we keep the balance condition
$$
\sum_{i=1}m_i>N-2.
$$
In one dimension it means $m_1>-1$ and then the theory done here is known to work, a detailed $L^1$ theory was developed in \cite{ERV88}. Note that we must write the equation as \ $u_t=(u_x/u^{1-m})_x$ \ to keep the singular parabolic type. A very curious case of non-uniqueness happens then.  In two dimensions we need  $m_1+m_2>0$ and $m_i<1$, which implies that $m_j>-1$. This lower bound holds in all dimensions $N\ge2$. We wonder if a large part of the above theory works. For the isotropic FDE see \cite{BBDGV} as a reference to this topic.

\smallskip

\noindent   {\bf Existence of solutions with locally finite data. } The local mass estimate of  Section \ref{local.mass}  allows one to develop a theory of locally bounded solutions for merely locally integrable data in $\ren $,   even exhibiting  ``burnt spots''. See our paper with Chasseigne \cite{ChV02} for  burnt spots, and also paper \cite{VaDarcy03} for receding pressure interfaces in 1D, which lead to solutions with arbitrarily expanding burnt spots.

\smallskip

\noindent    {\bf Functional topics.} Functional inequalities and entropy methods have been decisive in learning the fine details
of  FDE flows in the isotropic case.  There are plenty of earlier and current works, among them we will mention some references we are familiar with  and are related in some sense: \cite{BE84, BGL14, BBDGV, BBNS, CJM01, CV03, DK07, DS08, DPD02, DT13}. Nothing of the sort is known  in the anisotropic case yet.

Question: are there possible functional inequalities or entropies that reflect this anisotropy? Will some symmetrization work?

\smallskip

\noindent    {\bf Universal lower estimates.}  A priori bounds for positive solutions of the form
$$
\frac{\partial u}{\partial t}\ge - \frac{C u}{ t} \qquad \mbox{or} \qquad \frac{\partial u}{\partial t}\le \frac{C u}{ t}
$$
play a very useful role in the deep knowledge that has been gathered in the isotropic theory, \cite{BBDGV}. See our  estimate \eqref{asymp.V} for the VSS solutions. We wonder if there is any other related inequality in the anisotropic case for more general solutions. For a source reference to the original  estimates see \cite{BC81b}, see also \cite{VazSmooth}. In the other related models of nonlinear diffusion mentioned above much is known.

\smallskip

\noindent   {\bf Bounded domains. } We may study anisotropic diffusion in a domain of $\re^N$.
The local mass estimate of Section \ref{local.mass} allows to find VSS in bounded domains in the appropriate exponent range, but its role is to be ascertained. We do not know much about their properties or use.
A general theory of solutions in bounded domains is to be done. In the isotropic a very recent survey is \cite{BF24}.

\smallskip

 \noindent {\bf About the VSS.} Uniqueness of the VSS is not proved in this paper. We have uniqueness of the constructed VSS as a \sl minimal VSS solution\rm, which is enough for the theory of the paper. We guess that a {\sl maximal VSS} can also be shown to exist.

\smallskip

\noindent    {\bf Anisotropic gradient diffusion problems.}
A similar approach to the one displayed in this paper should work for the anisotropic $p$-Laplacian equation in the suitable fast diffusion range. The basic theory was done in \cite{FVV21} where the self-similarity was carefully analysed and related references can be found. General self-similar solutions are studied in \cite{Bive06} for the isotropic case.

Another interesting model of anisotropic diffusion uses the so-called doubly nonlinear  diffusion operator, mentioned and partially studied in  \cite{FVV21}.

\smallskip

\noindent {\bf Anisotropic fast diffusion with small exponents.} We have not studied what happens in the anisotropic FDE when $\overline{m}\le m_c$. The basic theory for the isotropic case can be found in \cite{VazSmooth}, see also \cite{BBDGV}. Conservation of mass does not hold and we may find extinction in finite time in many situations.

\smallskip

\noindent {\bf Solutions of both signs.} This interesting aspect is not touched.

\smallskip

\noindent {\bf Question.} Do we have explicit solutions in some cases?

\smallskip
%
%
%

\medskip


%
\nc

\section*{Acknowledgments}

The author was funded  by  grant PID2021-127105NB-I00 from MICINN (Spanish Government). He is an Honorary Professor at Univ. Complutense de Madrid and member of IMI. The paper was written during the ``2023 Thematic period on PDEs. Diffusion, Geometry, Probability and Free Boundaries'',
held at ICMAT-UAM, Madrid, Spain, in the period  June-December 2023.

\medskip



\begin{thebibliography}{999}


\bibitem{BE84} {\sc D. Bakry,   M. \'Emery.} Diffusions hypercontractives, Séminaire de probabilités,
XIX, 1983/84, Lecture Notes in Math., vol. 1123, Springer, Berlin, 1985, pp. 177--206.


\bibitem{BGL14} {\sc D. Bakry, I. Gentil, M. Ledoux.} Analysis and geometry of Markov diffusion operators. Grundlehren der mathematischen Wissenschaften, 348. Springer, Cham, 2014.

\bibitem{BC81b} {\sc Ph. B\'enilan, M.~G. Crandall.}
Regularizing effects of homogeneous evolution equations. {
Contributions to Analysis and Geometry}, (suppl.  to Amer. Jour.
Math.), Johns Hopkins Univ. Press, Baltimore, Md., 1981. Pp.
23-39.

\bibitem{Bive06} {\sc  M. F. Bidaut-V\'eron.} Self-similar solutions of the p-Laplace heat equation: the fast diffusion case, Pacific J. Math. 227 (2) (2006) 201--269.

\bibitem{BBDGV} {\sc  A. Blanchet, M. Bonforte, J. Dolbeault, G. Grillo, J. L. V\'azquez.} Asymptotics of the
fast diffusion equation via entropy estimates, Arch. Rational Mech. Anal. 191 (2009),
347--385.

%
\bibitem{BBNS} {\sc  M. Bonforte, J.~Dolbeault, B.~Nazaret, N. Simonov.} Stability in Gagliardo-Nirenberg-Sobolev inequalities. Flows, regularity and the entropy method. (171 pages). To appear in Memoirs AMS (2023).

\bibitem{BF24} {\sc  M. Bonforte, A. Figalli.}
The Cauchy-Dirichlet problem for the fast diffusion equation on bounded domains,
Nonlinear Analysis {\bf 239} (2024), 113394.



\bibitem{BonIagVaz2010} {\sc M. Bonforte, R. G. Iagar, J. L. V\'azquez. }
Local smoothing effects, positivity, and Harnack inequalities for the fast $p$-Laplacian equation.
Adv. Math. {\bf 224 } (2010), no. 5, 2151--2215.
\nc
%
\bibitem{BonSim19} {\sc M. Bonforte, N. Simonov.}  Quantitative a priori estimates for fast diffusion equations with Caffarelli-Kohn-Nirenberg weights. Harnack inequalities and H\"older continuity. Adv. Math. 345 (2019), 1075--1161.

\bibitem{BonSim23} {\sc M. Bonforte, N. Simonov.} Fine properties of solutions to the Cauchy problem for a fast diffusion equation with Caffarelli-Kohn-Nirenberg weights. Ann. Inst. H. Poincaré C Anal. Non Linéaire 40 (2023), no. 1, 1--59.

\bibitem{BonSimStan22} {\sc M. Bonforte, N. Simonov, D. Stan.} The Cauchy problem for the fast $p$-Laplacian evolution equation. Characterization of the global Harnack principle and fine asymptotic behaviour. J. Math. Pures Appl. (9) 163 (2022), 83--131.

\bibitem{BonVaz10} {\sc M. Bonforte, J. L. V\'azquez.}	Positivity, local smoothing, and Harnack inequalities for very fast diffusion equations, Adv. Math. 223 (2010), no. 2, 529--578.


\bibitem{BPT85} {\sc H. Brezis, L. A. Peletier, D. Terman.} A very singular solution of the heat equation
with absorption, Arch. Rational Mech. Anal. 96 (1985), 185--209.

\bibitem{CJM01} {\sc J. A. Carrillo, A. J\"ungel, P. A. Markowich, G. Toscani, A. Unterreiter.} Entropy dissipation
methods for degenerate parabolic problems and generalized Sobolev inequalities, Monatshefte f\"ur Mathematik, 133
(2001), pp. 1--82.

\bibitem{CV03} {\sc J. A. Carrillo, J. L. V\'azquez.} Fine asymptotics for fast diffusion equations, Communications in Partial
Differential Equations, 28 (2003), pp. 1023--1056.

\bibitem{ChV02} {\sc E. Chasseigne, J. L. Vazquez}.
Theory of Extended Solutions for Fast Diffusion Equations in
Optimal Classes of Data. Radiation from Singularities.
{ Archive Rat. Mech. Anal.} {\bf 164} (2002), 133--187.


\bibitem{DK07} {\sc P. Daskalopoulos, C. E. Kenig.}
Degenerate diffusions. Initial value problems and local regularity theory. EMS Tracts in Mathematics, 1. European Mathematical Society (EMS), Zürich, 2007.

\bibitem{DS08} {\sc P. Daskalopoulos, N. Sesum}. On the extinction profile of solutions to fast diffusion, Journal f\"ur die Reine
und Angewandte Mathematik. [Crelle's Journal], 622 (2008), pp. 95--119.

\bibitem{DPD02} {\sc M. Del Pino, J. Dolbeault.} Best constants for Gagliardo-Nirenberg inequalities and applications to non-linear diffusions, Journal de Math\'ematiques Pures et Appliqu\'ees. Neuvi\`eme S\'erie, 81 (2002), pp. 847--875.

\bibitem{DT13} {\sc J. Dolbeault, G. Toscani.} Improved interpolation inequalities, relative entropy and fast diffusion equations,
Annales de l'Institut Henri Poincar\'e. Analyse Non Lin\'eaire, 30 (2013), pp. 917--934.

\bibitem{ERV88} {\sc J. R. Esteban, A. Rodr\'{\i}guez, J. L. Vazquez.}
 A nonlinear heat equation with singular diffusivity.
{ Comm. Partial Diff. Eqs.} {\bf 13} (1988), 985--1039.

\bibitem{FVV21}  {\sc F. Feo,   J. L. V\'azquez, B. Volzone}. {
 Anisotropic $p$-Laplacian Evolution of Fast Diffusion type}.
 Adv. Nonlinear Stud. {\bf 21} (2021), no. 3, 523--555.

\bibitem{FVV23} {\sc F. Feo,   J. L. V\'azquez, B. Volzone}. {
Anisotropic Fast Diffusion Equations},
Nonlinear Anal. 233 (2023), Paper No. 113298, 43 pp.

\bibitem{H} {\sc E. Henriques.} {Concerning the regularity of the anisotropic porous medium equation}, J. Math. Anal. Appl. 377 (2011), n.2, 710--731.

\bibitem{HP} { \sc M. A. Herrero, M Pierre. } {The Cauchy problem for $u_t=\triangle u^m$ when $0<m<1$}, Trans. Amer. Math. Soc. 291 (1985), 145-158.

%

\bibitem{SJ05} {\sc B. H. Song, H.Y. Jian}. {Fundamental Solution of the Anisotropic Porous Medium Equation}, Acta Math. Sinica 21 n.5 (2005), 1183-1190.

\bibitem{SJ06} {\sc B. H. Song, H.Y. Jian}. {Solutions of the anisotropic porous medium equation in $\mathbb{R}^n$ under an $L^1$-initial value}, Nonlinear Analysis 64 (2006), 2098-2111.
%


\bibitem{Vascppme} {\sc J.~L.~V\'{a}zquez}.
{ Asymptotic behaviour for the Porous Medium Equation posed in the
whole space}. { J. Evolution Equations} {\bf 3} (2003),
67--118.

\bibitem{VaDarcy03} {\sc J.~L.~V\'{a}zquez}. Darcy's law and the theory of shrinking solutions of fast diffusion equations. SIAM J. Math. Anal. 35 (2003), no. 4, 1005--1028.

 \bibitem{Vlibro} {\sc J.~L.~V\'{a}zquez}. \emph{The Porous Medium Equation: Mathematical Theory},
 Oxford Mathematical Monographs. The Clarendon Press, Oxford University Press, Oxford  (2007).

\bibitem{VazSmooth} {\sc J.~L.~V\'azquez.} \emph{Smoothing and Decay Estimates for Nonlinear Diffusion Equations. Equations  of Porous Medium Type}. Oxford Lecture Series in Mathematics and Its Applications, vol. 33
(Oxford University Press, Oxford, 2006).

\bibitem{Vaz2014} {\sc J.~L.~V\'azquez,}
{Barenblatt  solutions and asymptotic behaviour for a  nonlinear fractional heat equation of porous medium type. }   J. Eur. Math. Soc. {\bf 16} (2014), 769--803.


\bibitem{Vaz20} {\sc J.~L.~V\'azquez.} {The evolution fractional p-Laplacian equation in $\ren$. Fundamental solution and asymptotic behaviour}. Nonlinear Anal. 199 (2020), 112034, 32 pp. \nc

\end{thebibliography}
\end{document}